\documentclass{amsart}

\newtheorem{theorem}{Theorem}[section]
\newtheorem*{maintheorem}{Theorem}
\newtheorem{lemma}[theorem]{Lemma}
\newtheorem{proposition}[theorem]{Proposition}
\newtheorem{corollary}[theorem]{Corollary}

\theoremstyle{definition}
\newtheorem{definition}[theorem]{Definition}
\newtheorem{example}[theorem]{Example}

\newtheorem{remark}[theorem]{Remark}
\newtheorem*{acknowledgement}{Acknowledgement}

\theoremstyle{remark}

\DeclareFontFamily{U}{wncy}{}
\DeclareFontShape{U}{wncy}{m}{n}{<->wncyr10}{}
\DeclareSymbolFont{mcy}{U}{wncy}{m}{n}
\DeclareMathSymbol{\Sh}{\mathord}{mcy}{"58}

\newcommand\mynote[1]{\marginpar{\ \\ \small \tt #1}}
\newcommand\comment[1]{\marginpar{\ \\ \small \sl \framebox{#1}}}
\newcommand\bel[1]{{\mynote{#1}}\begin{equation}\label{#1}}
\newcommand\mylabel[1]{\label{#1}}

\newcommand{\FF}{\mathbb{F}}

\newcommand{\PP}{\mathbb{P}}
\renewcommand{\AA}{\mathbb{A}}

\newcommand  {\shA}     {\mathcal{A}}

\newcommand  {\shC}     {\mathcal{C}}

\newcommand  {\shExt}   {\mathcal{E} \!\text{\textit{xt}}}

\newcommand  {\shF}     {\mathcal{F}}
\newcommand  {\shG}     {\mathcal{G}}

\newcommand  {\shHom}   {\mathcal{H}\!\text{\textit{om}}}
\newcommand  {\shI}     {\mathcal{I}}
\newcommand  {\shJ}     {\mathcal{J}}

\newcommand  {\shM}     {\mathcal{M}}

\newcommand  {\shL}     {\mathcal{L}}

\newcommand  {\shT}     {\mathcal{T}}

\newcommand  {\ida}      {\mathfrak{a}}

\newcommand  {\Cl}      {\operatorname{Cl}}

\renewcommand{\cong}    {\equiv}

\newcommand  {\depth}   {\operatorname{depth}}

\newcommand  {\edim}    {\operatorname{edim}}

\newcommand  {\Et}      {{\operatorname{Et}}}

\newcommand  {\Gal}     {\operatorname{Gal}}

\newcommand  {\Hom}     {\operatorname{Hom}}

\renewcommand  {\ker }  {\operatorname{kern}}

\newcommand  {\loc}     {{\operatorname{loc}}}
\newcommand  {\dirlim}  {\varinjlim}

\newcommand  {\lra}     {\longrightarrow}

\newcommand  {\maxid}   {\mathfrak{m}}

\renewcommand{\O}       {\mathcal{O}}

\newcommand  {\Pic}     {\operatorname{Pic}}

\newcommand  {\pd}      {{\operatorname{pd}}}

\newcommand  {\quadand} {\quad\text{and}\quad}

\newcommand  {\ra}      {\rightarrow}

\newcommand  {\red}     {{\operatorname{red}}}

\newcommand  {\Set}     {(\text{Set})}

\newcommand  {\sh}      {{\operatorname{sh}}}

\newcommand  {\SL}      {\operatorname{SL}}

\newcommand  {\Spec}    {\operatorname{Spec}}

\newcommand  {\trdeg}   {\operatorname{trdeg}}

\def\mydate{\number\day\space\ifcase\month \or January\or February\or March\or 
April\or May\or June\or July\or
August\or September\or October\or November\or December\fi \space\number\year}

\DeclareFontFamily{U}{wncy}{}
\DeclareFontShape{U}{wncy}{m}{n}{<->wncyr10}{}
\DeclareSymbolFont{mcy}{U}{wncy}{m}{n}
\DeclareMathSymbol{\Sh}{\mathord}{mcy}{"58}

\usepackage{amscd,amssymb,amsmath}

\begin{document}

\title[Singularities appearing on   generic fibers]
      {Singularities appearing on   generic fibers
of   morphisms between smooth schemes}

\author[Stefan Schroer]{Stefan Schr\"oer}
\address{Mathematisches Institut, Heinrich-Heine-Universit\"at,
40225 D\"usseldorf, Germany}
\curraddr{}
\email{schroeer@math.uni-duesseldorf.de}


\dedicatory{1 August 2006}

\begin{abstract}
I give  various criteria for singularities to appear on geometric
generic fibers of morphism between smooth schemes in positive characteristics.
This involves   local fundamental groups,
jacobian ideals, projective dimension, tangent and cotangent sheaves, and 
the effect of Frobenius.
As an application, I determine which rational double points do 
appear on geometric generic fibers.
\end{abstract}

\maketitle
\tableofcontents

\section*{Introduction}

The goal of this paper is to explore the   structure of singularities
that   occur on generic fibers in positive characteristics. As an application of our
general results,
we shall determine which rational double points do and which do not occur on generic fibers.

In some sense, the starting point   is Sard's Lemma from differential topology.
It states that the   critical values of a differential map between
differential manifolds form a set of measure zero. As a consequence, any general fiber of   a differential
map is itself a differential manifold. The analogy  in algebraic geometry
is as follows:  
Let $k$ be an algebraically closed ground field of characteristic $p\geq 0$,
and suppose $f:S\ra B$ is a morphism   between smooth integral schemes.
Then the generic fiber $S_\eta$ is a regular   scheme of finite type over the function field
$E=\kappa(\eta)$.

In characteristic zero, this implies that $S_\eta$ is smooth over $E$.
Moreover, the absolute Galois group $G=\Gal(\bar{E}/E)$ 
acts on the geometric generic
fiber $S_{\bar{\eta}}$ with quotient isomorphic to $S_\eta$.
In other words, to understand the generic fiber, it suffices to understand
the geometric generic fiber, which is again smooth over an algebraically closed field,
together with its Galois action.

The situation is more complicated in characteristic $p>0$.
The reason is that over nonperfect fields, the notion of regularity, which depends only on the scheme and not
on the structure morphism, 
is weaker that the notion of geometric regularity, which coincides with formal smoothness.
Here it easily happens that the geometric generic fiber $S_{\bar{\eta}}$
acquires singularities. This special effect already plays a crucial role in
the extension of Enriques classification of surface to positive characteristics: 
As Bombieri and Mumford \cite{Bombieri; Mumford 1976} showed, 
there are \emph{quasielliptic fibrations} for $p=2$ and $p=3$, 
which are analogous to elliptic fibrations, but have a cusp   on the geometric generic fiber.

We call a proper morphism  $f:S\ra B$ of smooth algebraic schemes a \emph{quasifibration}
if $\O_B=f_*(\O_S)$, and the generic fiber $S_\eta$ is not smooth.
The existence of quasifibrations should by no means be viewed as pathological.
Rather, they involve some fascinating geometry, and apparently  offer  new
 freedom to achieve   geometrical constructions
that are impossible in characteristic zero. 
The theory of quasifibrations, however, is still in its infancy.
In \cite{Kollar 1991}, Remark 1.2, Koll\'ar asks whether or not there are Fano contractions
on threefolds whose geometric generic fibers are nonnormal del Pezzo surfaces.
Some results in this direction appear in \cite{Schroeer 2005}.
Mori and Saito \cite{Mori; Saito 2003} studied Fano contractions
whose geometric generic fibers are nonreduced quadrics.
Examples of quasifibrations involving minimally elliptic singularities
appear in \cite{Schroeer 2006}, in connection with Beauville's generalized
Kummer varieties.

Of course, nonsmoothness of 
the generic fiber $S_\eta$ leads to unusual  complications.
However, singularities appearing on the geometric generic fiber
$S_{\bar{\eta}}$ are not arbitrary. First and formost, they are  locally of 
complete intersection, whence many powerful methods from commutative algebra apply.
However,  they satisfy far more restrictive
conditions, and the goal of this paper is to analyze these conditions.
Hirokado \cite{Hirokado 2004} started such an analysis, and characteriszed
those rational double points in odd characteristics that appear
on geometric generic fibers. His approach was to study the closed fibers
$S_b$, $b\in B$ and their deformation theory. Our approach is  somewhat different:
We look at the generic fiber $S_\eta$ and deliberately  work
over the function field $\kappa(\eta)$.

In fact, we will mainly work in the following abstract setting:
Given an field $F$ in characteristic $p>0$,
and a subfield $E$ so that the field extension $E\subset F$ is purely inseparable,
we consider $F$-schemes $X$ of finite type that descend to regular $E$-schemes $Y$, that
is $X\simeq Y\otimes_E F$. Our first results on such schemes $X$ are as follows:
In codimension two, the local fundamental groups are trivial, and the torsion of the local
class groups are $p$-groups. Moreover, the Tjurina numbers are divisible by $p$, the stalks of the jacobian
ideal have finite projective dimension, and the tangent sheaf $\Theta_X$ is locally free in codimension two.
These conditions are comparatively straightforward, but already give strong conditions on the
singularities. The following restriction on the cotangent sheaf was a bit
of a surprise to me, and is  the first main result of this paper:

\begin{maintheorem}
If an $F$-scheme $X$ descends to a regular scheme, then
the for each point $x\in X$ of codimension two, the stalk $\Omega^1_{X/F,x}$ contains an invertible direct summand.
In other words, $\Omega^1_{X/F,x}\simeq\O_{X,x}\oplus M$ for some $\O_{X,x}$-module $M$.
\end{maintheorem}

Note that any torsion free module of finite type over an integral local noetherian ring is an extension of some ideal
by a free modul. Such extensions, however, do not necessarily split, so the preceding result puts
a nontrivial condition on the cotangent sheaf.

As an application of all these results, we shall determine which rational double points
appear on surfaces  descending to regular schemes, and which do not.
This was already settled by Hirokado \cite{Hirokado 2004} in the case of odd characteristics.
Recall that Artin \cite{Artin 1977} classified    rational double points in positive characteristics.
Here the isomorphism class in not merely determined by a Dynkin diagram, but sometimes on  some additional integral parameter $r$ 
(this is   the case for $p=2,3,5$). It turns out that the situation is most
challenging in characteristic two:
Besides the $A_n$-singularities, which behave as in characteristic zero, there are the following isomorphism classes:
$$
D_n^r,\quad \text{with $ 0\leq r\leq \lfloor n/2\rfloor -1$}, \quadand E_6^0,E_6^1,E_7^0,\ldots,E_7^3,E_8^0,\ldots, E_8^4.
$$
For simplicity, 
I state the second main result of this paper only in characteristic two,
which answers a question of Hirokado \cite{Hirokado 2004}:

\begin{maintheorem}
In characteristic two, the rational double points that apprear on surfaces  descend  to 
regular surfaces  are the following: $A_{2^e-1}$ with $e\geq 1$, and $D_n^0$ with $n\geq 4$,
and $E_n^0$ with $n=6,7,8$.
\end{maintheorem}

Note that this   includes, but does not coincide with,
all rational double points that are purely inseparable
double coverings of a smooth scheme. 
Shepherd--Barron \cite{Shepherd-Barron 1996} calls them  \emph{special},
and showed that they play an important role in the geometry of surfaces in characteristic two. 
The $D_{2m+1}^0$-singularities are not special but nevertheless descend to regular surfaces. 
I also want to point out that all members of our list
have a tangent sheaf that is locally free, but there are   other rational
double points with locally free tangent sheaf.

Here is an outline of the paper:
In Section \ref{purely inseparable descend}, I set up notation
and give some elementary examples and results.
In the next four Sections, we analyze $F$-schemes $X$ that descend
to regular $E$-schemes $Y$.
In Section \ref{local fundamental group} we treat the local fundamental
groups.
In Section \ref{residue fields divisibility}, we shall see that integer-valued invariants attached
to the singularities of $X$ like Tyurina numbers are multiples of $p$.
Section \ref{finite projective dimension} deals with finite projective dimension
of sheaves obtained from the cotangent sheaf $\Omega^1_{X/F}$. Here we will also see that the tangent sheaf $\Theta_X$ is locally free in codimenson two.
Our first main result appears in Section \ref{invertible summands cotangent}: We show that the stalks
in codimension two of the cotangent sheaf $\Omega_{X/F}$ contains invertible summand.
This gives a powerful and easy-to-use criterion for complete intersections
defined by a single equation.
We will exploit this in Section \ref{rational double points}: Here
we determine which rational double points are possible on schemes
descending to regular schemes.

\begin{acknowledgement}
I wish to thank Holger Brenner, Torsten Ekedahl, Almar Kaid, Manfred Lehn, and my wife Le Van  
for helpful discussion. I have used the computer algebra system Magma for some computations
of lengths of Artin rings.
\end{acknowledgement}

\section{Purely inseparable descend}
\mylabel{purely inseparable descend}

Let $F$ be a  field of characteristic $p\geq 0$, and $E\subset F$ be a subfield.
For each $E$-scheme $Y$, base-change gives an $F$-scheme $X=Y\otimes_E F$.
Conversely, given an $F$-scheme $X$, one may ask whether or not there is 
an $E$-scheme $Y$ with the property $X\simeq Y\otimes_E F$.
If this is the case, we shall say that $X$ \emph{descends} along $E\subset F$.
Then the fiber product $R=X\times_Y X$ defines an equivalence relation on $X$,
and we may view $Y=X/R$ as the coresponding quotient.
The topic of this paper are schemes that are not regular itself but descend to regular schemes.
A bit of terminology:

\begin{definition}
We say that a locally noetherian $F$-scheme $X$  \emph{descends to a 
regular scheme} if there is a subfield $E\subset F$ and a regular $E$-scheme $Y$
with $X\simeq Y\otimes_E F$.
\end{definition}

This notion is of little interest for us in characteristic zero: Then any locally noetherian scheme   descending
to a regular scheme is itself regular, which follows from \cite{EGA IVb}, Proposition 6.7.4.
Therefore we assume from now on that we are in characteristic $p>0$.
Throughout,  $X$ ususally denotes a locally noetherian $F$-scheme.

\begin{lemma}
\mylabel{purely inseparable}
Suppose that  the $F$-scheme $X$ descends to a regular scheme. Then there is a subfield $E$ inside $F$
and a regular $E$-scheme $Y$ with $X\simeq Y\otimes_E F$ so that the  field extension $E\subset F$ is purely inseparable.
\end{lemma}

\proof
We have $X\simeq Y_0\otimes_{E_0} F$ with a regular $E_0$-scheme $Y_0$ for some
subfield $E_0\subset F$. Choose a transcendence basis $f_\alpha\in F$, $\alpha\in I$ for this field extension,
and let $E$ be the separable algebraic closure of $E_0(f_\alpha)_{\alpha\in I}$ inside $F$.
Then the field extension $E\subset F$ is purely inseparable, and the
scheme $Y=Y_0\otimes_{E_0} E$ remains regular, by \cite{EGA IVb}, Proposition 6.7.4.
\qed

\medskip
Let me discuss an example to see how this might happen.
Suppose $v_0,v_1,\ldots,v_n$ is  a collection of indeterminates, and let $f\in E[v_0,v_1,\ldots,v_n]$
be a polynomial of the form $f=v_0^q+g$, where $q=p^e$ is a power of the characteristic,
and $g\in (v_1,\ldots,v_2)^2$ does not involve the indeterminate $v_0$ and has neither constant nor linear terms.
Suppose there is an element $\lambda\in E$    
that is not a  $p$-th power in $E$, but even becomes a $q$-th power in $F$.
Consider the two affine schemes
$$
Y=\Spec E[v_0,v_1,\ldots,v_n]/(f-\lambda)\quadand
X=\Spec F[v_0,v_1,\ldots,v_n]/(f).
$$
Let $y\in Y$ be the  closed point  defined by $v_1=\ldots=v_n=0$. Note
that this is not a rational point. Rather, its residue field is $\kappa(y)=E(\lambda^{1/q})$.
Let $x\in X$ be the rational point given by $v_0=v_1=\ldots=v_n=0$.

\begin{proposition}
\mylabel{regular singular}
With the preceding assumptions, the local ring $\O_{Y,y}$ is regular,    the local ring $\O_{X,x}$ is not regular,
and $\O_{X,x}\simeq \O_{Y,y}\otimes_EF$ holds. Hence some open neighborhood of the singularity $x\in X$
descends to an regular scheme.
\end{proposition}

\proof
To check that $\O_{Y,y}$ is regular,
set $\AA^n_E:=\Spec E[v_1,\ldots,v_n]$  and consider the canonical morphism $\varphi:Y\ra\AA^n_E$.
This is flat of degree $q$, and the fiber over the origin $0\in\AA^n_E$ is isomorphic
to the spectrum of $E[v_0]/(v_0^q-\lambda)$, which is a field.
In particular, the base $\AA^n_E$ and the fiber $\varphi^{-1}(0)$ are both regular schemes.
According to \cite{EGA IVa}, Proposition 17.3.3, this implies that the local ring $\O_{Y,y}$ is regular. 

By our assumptions on $f$, we  have 
$$
x\in\Spec k[v_0,v_1,\ldots,v_n]/(v_0,v_1,\ldots,v_n)^2\subset X.
$$
Whence $\edim(\O_{X,x})\geq n+1>n=\dim(\O_{X,x})$. It follows that the local ring $\O_{X,x}$ is
not regular. Finally, the  automorphism of $F[v_0,v_1,\ldots,v_n]$ defined by the substitution $v_0\mapsto v_0+\lambda^{1/q}$ maps $f-\lambda$
to $f$, whence $X\simeq Y\otimes_E F$, and this isomorphism gives  $\O_{X,x}\simeq \O_{Y,y}\otimes_EF$.
\qed

\begin{remark}
 There are singularities of   different structure that descend to regular schemes,
for example   rational double points of type $D_{2m+1}^0$ in characteristic two.
Compare Section \ref{rational double points}.
\end{remark}

Schemes descending to regular schemes are not arbitrary. Rather, they satisfy several strong conditions.
Recall that a locally noetherian scheme $S$ is called \emph{locally of complete intersection} if 
for each point $s\in S$, the formal completion   is of the form
$\O_{S,s}^\wedge\simeq R/(f_1,\ldots,f_r)$,
where $R$ is a regular complete local noetherian ring, and $f_1,\ldots,f_r\in R$ is a regular sequence.

\begin{proposition}
\mylabel{complete intersection}
Suppose a locally noetherian $F$-scheme $X$ descends to regular scheme. Then $X$ is locally of complete intersection.
\end{proposition}

\proof
Regular schemes are obviously locally of complete intersection.
According to \cite{EGA IVd},  Corollary 19.3.4,
the property of being locally of complete intersection
is preserved under base field extensions.
\qed

\medskip
Suppose $k$ is a ground field of characteristic $p>0$,
and let $S$ be a smooth   $k$-scheme. Suppose $f:S\ra B$ is a morphism
onto an integral $k$-scheme of finite type, and  $\eta\in B$ be the generic point.
Then the generic fiber $S_\eta$ is regular, because the morphism $S_\eta \ra S$
is flat with regular fibers. It follows that the geometric generic fiber $S_{\bar{\eta}}$, which is of finite type over $F=\overline{\kappa(\eta)}$,
descends to a regular scheme.
Conversely:

\begin{proposition}
 Let $X$ be an connected $F$-scheme of finite type that descends to a regular scheme.
Then there are smooth connected $\FF_p$-schemes of finite type $S$ and $B$,
a dominant morphism $f:S\ra B$, and an inclusion of fields $\kappa(\eta)\subset F$
with $X\simeq S_\eta\otimes_{\kappa(\eta)} F$. Here $\eta\in B$ denotes the generic point.
\end{proposition}

\proof
Write $X=Y\otimes_E F$ for some regular $E$-scheme $Y$. The $E$-scheme $Y$ is of
finite type, because the $F$-scheme $X$ is (\cite{SGA 1}, Expos\'e VIII, Proposition 3.3).
We therefore may assume that the field $E$ is finitely generated over its
prime field $\FF_p$, according to \cite{EGA IVc}, Theorem 8.8.2. 
Whence $E=\kappa(B)$ is the function field of some
integral $\FF_p$-scheme $B$ of finite type. Moreover,   $Y$ is isomorphic to  the
generic fiber of some dominant morphism $f:S\ra B$. Shrinking $S$ and $B$,
we may assume that $S$  and $B$ are regular. Then they are even smooth, because the
finite field $\FF_p$ is perfect.
\qed

\medskip
Now suppose that $S$ is smooth and proper.
We call a morphism $f:S\ra B$ onto another proper $k$-scheme $B$ a \emph{quasifibration}
if $\O_B=f_*(\O_S)$, and the generic fiber $S_\eta$ is not smooth.
Whence the geometric generic fiber $S_{\bar{\eta}}$ is singular,
but descends to a regular scheme.
The most prominent example for quasifibrations are the quasielliptic surfaces.

\section{Local fundamental groups}
\mylabel{local fundamental group}

Let $F$ be a field of characteristic $p>0$, and $X$   a locally noetherian $F$-scheme.
To avoid endless repetition of the same hypothesis,  \emph{we suppose throughout this section  that  $X$ descends to a regular scheme}.
In other words, there is a subfield $E$ so that this field extension $E\subset F$
is purely inseparable, and a regular $E$-scheme $Y$ with $X=Y\otimes_E F$. 
 The goal of this section
is to find out what this implies for the fundamental groups attached to $X$. The main idea
is to use  the Zariski--Nagata Purity
Theorem.

To start with, let me recall the definition of
fundamental groups in algebraic geometry. 
Suppose that $S$ is connected scheme. We shall denote $\Et(S)$ the category of finite \'etale morphism $S'\ra S$.
This is a Galois category, hence equivalent to the category $\shC(\pi)$ of finite sets
endowed with a continuous action of a profinite group $\pi$ as follows:
Suppose $a:\Spec(\Omega)\ra S$ is a base point, for some separably closed field $\Omega$.
This yields a fiber functor $\Et(S)\ra\Set$, $S'\mapsto S'_a$.
The \emph{fundamental group} $\pi_1(S,a)$ is defined as the automorphism group of this fiber
functor. As explained in \cite{SGA 1},  Expos\'e V, Section 4, this makes the fundamental group profinite,
and the fiber functor yields an equivalence $\Et(S)\ra\shC(\pi_1(S,a))$.
The choice of different base points yields isomorphic fundamental groups, but the isomorphism
is only unique up to inner automorphisms. By abuse of notation, we sometimes write
$\pi_1(S)$ when  speaking about group-theoretical properties that depend only on   isomorphism
classes  of  groups.

Now back to our situation $X=Y\otimes_E F$. 
Let $a$ be a geometric point on $X$, and $b$ the induced geometric point on $Y$.
We have the following basic fact:

\begin{lemma}
\mylabel{fundamental group}
The canonical homomorphism $\pi_1(X,a)\ra\pi_1(Y,b)$ is an isomorphism of topological groups.
\end{lemma}

\proof
By assumption, the field extension $E\subset F$ is purely inseparable.
Consider first the special case that $E\subset F$ is finite as well.
Then \cite{SGA 1}, Expos\'e IX, Theorem 4.10 tells us that the pull-back   functor
$\Et(Y)\ra\Et(X)$, $Y'\mapsto Y'\times_Y X$ is an  equivalence of categories.
In the general case write $F=\bigcup F_\lambda$ as a union of subfields that are finite over $E$.
Using \cite{EGA IVc}, Theorem 8.8.2, one infers that the pull-back functor $Y'\mapsto Y'\times_Y X$
is still an equivalence of categories. The assertion follows.
\qed

\medskip
This leads   to the following result:

\begin{theorem}
\mylabel{zariski nagata}
Suppose our $F$-scheme $X$ is local and henselian, with closed point $x\in X$. Let $A\subset X$ be a closed
subset of codimension $\geq 2$, and $U\subset X$ its complement.
Then $U$ is connected, and its fundamental group $\pi_1(U)$ is  isomorphic to the absolute Galois group of the residue field $\kappa(x)$.
\end{theorem}

\proof
If follows from Proposition \ref{complete intersection} that the scheme $X$ is Cohen--Macaulay, 
hence $\depth(\O_{X,a})\geq 2$ for all points $a\in A$. By Hartshorne's Connectedness Theorem (\cite{SGA 2}, Expos\'e III, Theorem 3.6), the complement $U$ is connected.
Let $V\subset Y$ be the open subset corresponding to $U\subset X$.
Consider the following chain of restriction and pull-back functors
$$
\Et(U)\stackrel{r_1}{\longleftarrow} \Et(V) \stackrel{r_2}{\longleftarrow} \Et(Y)
\stackrel{r_3}{\lra}\Et(X) \stackrel{r_4}{\lra}\Et(x).
$$
According to Lemma \ref{fundamental group}, the pull-back
functors $r_1,r_3$ are equivalences. By assumption, $Y$ is regular, so we may 
apply the Zariski--Nagata Theorem (see \cite{SGA 2}, Expos\'e X, Theorem 3.4 for a scheme-theoretic proof) and deduce 
that  the restriction functor $r_2$ is an equivalence as well.
Since $X$ is henselian, the restriction functor $r_4$ is an equivalence,
as explained in \cite{EGA IVd}, Proposition 18.5.15.
The assertion follows.
\qed

\medskip
Recall that a local scheme is called \emph{strictly henselian} if
it is henselian, and the residue field of the closed point is separably closed.
Hence we have:

\begin{corollary}
\mylabel{strictly henselian}
Assumption as in Theorem \ref{zariski nagata}. If $X$ is even strictly henselian,
then $U$ is connected and simply connected.
\end{corollary}

\medskip
Recall that for a local scheme $S$, the \emph{local fundamental group} $\pi_1^\loc(S,a)$
is defined as the fundmental group of the complement of the closed point $s\in S$,
where $a$ denotes a base point on the scheme $S\setminus\left\{s\right\}$.
We may apply the preceding corrollary to the strict henselization $\Spec(\O_{X,x}^\sh)$
of a point $x\in X$, where $X$ is an arbitrary locally noetherian scheme.

\begin{corollary}
\mylabel{local fundamental}
Let $x\in X$ be a point of codimension $\geq 2$.
Then $\pi_1^\loc(\O_{X,x}^\sh)=0$.
\end{corollary}

\proof
In order to apply Corollary \ref{strictly henselian}, we merely have to verify that the henselization $S:=\Spec(\O_{X,x}^\sh)$ descends to a regular scheme.
Indeed: Let $y\in Y$ be the image of $x\in X$. The canonical map $\O_{Y,y}^\sh\ra\O_{X,x}^\sh$ coming from the universal property of strict henselizations  on $Y$ induces a map $\varphi:\O_{Y,y}^\sh\otimes_E F\ra\O_{X,x}^\sh$.
The scheme  $\O_{Y,y}^\sh\otimes_E F$ is clearly local with separably closed residue field, and it follows
from \cite{EGA IVd}, Theorem 18.5.11 that it is   henselian as well. Whence the universal
property of strict henselization on $X$ yields the desired inverse map
$\psi:\O_{X,x}^\sh\ra \O_{Y,y}^\sh\otimes_E F$. 
\qed

\medskip
For an arbitrary locally noetherian scheme $S$, we define the \emph{class group} as $\Cl(S)=\dirlim \Pic(U)$, where
the filtered direct limit runs over all open subsets $U\subset S$ so that
$\depth(\O_{S,s})\geq 2$ for all $s\in S\setminus U$.
In case $S$ is normal, this is indeed the group of Weil divisors modulo principal divisors
(\cite{SGA 2}, Expos\'e XI, Proposition 3.7.1). For  nonnormal or even nonreduced  schemes, however, the above definition seems to be more suitable.

\begin{corollary}
\mylabel{class group}
Let $x\in X$ be a point of codimension $\geq 2$.
Then the torsion part of the class group $\Cl(\O_{X,x}^\sh)$ is a $p$-group.
\end{corollary}

\proof
Set $S=\Spec(\O_{X,x}^\sh)$. It follows from \ref{complete intersection} that
the scheme $S$ is Cohen--Macaulay. Whence the points $s\in S$ with $\depth(\O_{S,s})\geq 2$
are precisely the points of codimension $\geq 2$. Seeking a contradiction, we suppose
that there is an open subset $U\subset S$ whose complement has codimension $\geq 2$,
and an invertible $\O_U$-module $\shL$ whose order $l> 1$ in $\Pic(U)$ is finite but  not a $p$-power.
Passing to  a suitable multiple, we may assume that $l$ is a prime number different from $p$.
 Choosing a trivialization of $\shL^{\otimes l}$, we get an $\O_U$-algebra structure on
the coherent $\O_U$-module
$\shA=\O_U\oplus\shL\oplus\ldots\oplus\shL^{\otimes (l-1)}$. Its relative spectrum $U'=\Spec(\shA)$
is a finite \'etale Galois covering $U'\ra U$ of degree $l$.
This covering is nontrivial.
Whence $U$ is not simply connected. On the other hand, Corollary 
\ref{strictly henselian} tells us that $U$ is simply connected, contradiction.
\qed

\medskip
We now specialize the preceding result to the case of normal surface singularities.
Let $x\in X$ be a point of codimension two, and assume that the local noetherian ring $\O_{X,x}$ is  normal.
Let  $r:\tilde{X}\ra X$ be  a resolution of singularities (whose existence we tacitly assume), and 
$D=r^{-1}(x)_\red$ the reduced exceptional divisor.
Consider the connected component $\Pic^0_{D/\kappa(x)}$ of the Picard scheme, which parameterizes  invertible sheaves
on $D$ with zero degree on each integral component.
We call the singularity $x\in X$ \emph{unipotent} if the Picard scheme $\Pic^0_{D/\kappa(x)}$ is unipotent.
It is easy to verify that this does not depend on the choice of resolution of singularities.
Note also that to check this,  we may  use  any divisor $D\subset \tilde{X}$ whose support is the full fiber.
Rational singularities are examples of  unipotent singularities.

\begin{corollary}
\mylabel{unipotent singularity}
Let $x\in X$ be a point of codimension two.
Then the singularity $x\in X$
is unipotent.
\end{corollary}

\proof
Seeking a contradiction, we assume that the singularity $x\in X$ is not unipotent, such that
the group scheme $\Pic^0_{D/k}$ is not unipotent, where $k=\kappa(x)$. Choose a separable
closure $k\subset k^s$, and set $D^s=D\otimes_k k^s$.
There must be an invertible sheaf $\shL_0\in\Pic(D^s)$ of finite order $l$ prime to $p$.
It extends to an invertible  sheaf $\shL$ on the strict henselization $\O_{X,x}^\sh$ 
of order $l$, as explained in \cite{Schroeer 2001}, Lemma 2.2.
Restricting to the complement $U\subset \Spec(\O_{X,x}^\sh)$ of the closed point, we obtain an invertible sheaf  $\shL_U\in\Pic(U)$
of order $l$, in contradiction to Corollary \ref{class group}.
\qed

\begin{corollary}
\mylabel{nonsmooth components}
Assumptions as in the preceding corollary.
Let $D\subset\tilde{X}$ be the reduced fiber over the singularity $x\in X$,
and $D_i\subset D$ be its integral components. Then the schemes $D_i$ are geometrically unibranch,
and their intersection graph  is a tree. If $D_i$ is smooth, then it
is isomorphic to a quadric in $\PP^2$, and  $h^1(\O_{D_i})=0$.
\end{corollary}

\proof
If one of the conditions would not hold, the Picard scheme $\Pic_{D/\kappa(x)}$ would 
not be unipotent, as explained in \cite{Bosch; Luetkebohmert; Raynaud 1990}, Chapter 9.
\qed

\begin{remark}
Recall that a normal point $x\in X$ of codimension two is called a \emph{simple elliptic singularity}
if the reduced fiber $D\subset\tilde{X}$ is an elliptic curve.
The preceding corollary tells us that such singularities do not descend to   $E$-schemes that are regular.
Hirokado \cite{Hirokado 2004} showed   this by using the equations for  simple elliptic singularities.
\end{remark}

\section{Residue fields and $p$-divisibility}
\mylabel{residue fields divisibility}

In this section  $X$ denotes an  $F$-schemes of finite type.
We continue to assume throughout that $X$ descends to a regular scheme.
In other words, there is a subfield $E $ so that this field extension $E\subset F$ is
purely inseparable,  and a regular $E$-scheme $Y$ with $X=Y\otimes_EF$. 

Let $x\in X$ be a   point, and $y\in Y$ be its image.  By our assumption, the local ring $\O_{Y,y}$ is regular.
What can be said about the local ring $\O_{X,x}$ and the inclusion $\O_{Y,y}\subset\O_{X,x}$?
We start with an observation on the residue fields. Consider the following
commutative diagram
$$
\begin{CD}
\kappa(y) @>>> \kappa(x)\\
@AAA @AAA\\
E @>>> F
\end{CD}
$$
of fields, and choose an algebraic closure $\kappa(x)\subset\Omega$.

\begin{proposition}
\mylabel{linearly disjoint}
Suppose $x\in X$ is not regular. Then the subfields $\kappa(y),F\subset\Omega$
are not linearly disjoint over $E$.
\end{proposition}

\proof
Seeking a contradiction, we assume that the subfields $\kappa(y),F\subset\Omega$
are linearly disjoint. By definition, this means that the canonical map $\kappa(y)\otimes_E F\ra\Omega$
is injective, such that the ring $R=\kappa(y)\otimes_E F$ is integral.
Using that the ring extension $\kappa(y)\subset R$ is integral, we infer that the   ring $R$ is actually a field.
It follows that the fiber $X_y=\Spec(R)$ of the projection $X\ra Y$ is a regular scheme.
By assumptions, $y\in Y$ is regular. By \cite{EGA IVb}, Corollary 6.5.2,
this implies that $x\in X$ is regular, contradiction.
\qed

\medskip
This has the following numerical consequence:

\begin{proposition}
\mylabel{not separable}
Suppose $x\in X$ is not regular. Then the   field extension $E\subset\kappa(y)$ is not separable.
If in addition $x\in X$ is closed, then the characteristic $p$ divides the
degree $[\kappa(y):E]$ of the image point $y\in Y$.
\end{proposition}

\proof
Suppose that $E\subset \kappa(y)$ is separable. Then the ring $R=\kappa(y)\otimes_E F$
is reduced. It is also irreducible, because $E\subset F$ is purely inseparable.
It follows that $R$ is integral,   whence the canonical map $R\ra \Omega$ must be injective,
and $\kappa(y),F\subset\Omega$ are linearly disjoint, contradiction.

If the point $x\in X$ is closed, then the field extension $\kappa(y)\subset E$
is finite, and its degree is of the form $[\kappa(y):E]=p^e[\kappa(y):E]_s$,
where the second factor is the  separable degree, and $e\geq 0$ is an integer.
Since $E\subset\kappa(y)$ is not separable, we must have $e>0$. 
\qed

\medskip
We shall apply this result frequently in the following form:

\begin{lemma}
\mylabel{multiple p}
Let $\shF_Y$ be a coherent zero-dimensional sheaf on $Y$ supported 
on the nonsmooth locus,
and $\shF_X=\shF_Y\otimes_EF$ the corresponding coherent sheaf on $X$.
Then the number $h^0(X,\shF_X)$ is a multiple of $p$.
\end{lemma}

\proof
Recall that $h^0(X,\shF_X)$ is the dimension of   $H^0(X,\shF_F)$ viewed as an
$F$-vector space.
By flat base-change, we have $h^0(Y,\shF_Y)=h^0(X,\shF)$.
The coherent zero-dimensional sheaf $\shF_Y$ has a Jordan--H\"older series with
factors   isomorphic to residue fields $\kappa(y)$, where
$y\in Y$ are nonsmooth closed points. By Proposition \ref{not separable},
the dimension of the $E$-vector space $\kappa(y)$ is a multiple
of $p$, and the assertion follows.
\qed

\medskip
There are many such coherent sheaves  $\shF_X=\shF_Y\otimes_E F$ coming from the sheaf
of K\"ahler differentials, provided that
the nonsmooth locus is zero-dimensional. Let me mention a few:

\begin{proposition}
\mylabel{exterior powers}
Suppose the nonsmooth locus of $X$ is zero-dimensional.
Then for all $m>\dim(X)$, the numbers $h^0(X,\Omega^m_{X/k})$ are multiplies of $p$.
\end{proposition}

\proof
The coherent sheaf $\Omega^1_{Y/k}$ is locally free of rank $d=\dim(Y)$
on the smooth locus. For $m>\dim(Y)$,   the coherent sheaf 
$\shF_Y=\Omega^m_{Y/k}$   is therefore supported on
the nonsmooth locus. Hence Lemma \ref{multiple p} applies.
\qed

\medskip
Next, consider the coherent sheafs $\shT_m,\shT'_m$ defined by the exact sequence
$$
0\lra \shT_m\lra\Omega^m_{X/k}\lra(\Omega^m_{X/k})^{\vee\vee}\lra\shT'_m\lra 0,
$$
where the map in the middle is the canonical evaluation map into the bidual.

\begin{proposition}
 \mylabel{kernel cokernel}
Suppose the nonsmooth locus of $X$ is zero-dimensional.
Then for all $m\geq 0$, the numbers $h^0(X,\shT_m)$ and $h^0(X,\shT_m')$
are multiples of $p$.
\end{proposition}

\proof
Follows as above from Lemma \ref{multiple p}.
\qed

\medskip
Recall that the \emph{jacobian ideal} is defined as the
$d$-th Fitting ideal of $\Omega^1_{X/F}$, where $d=\dim(X)$. Here we tacitely assume
that $X$ is equidimensional. We call the closed subscheme $X'\subset X$ defined by
the jacobian ideal the \emph{jacobian subscheme}. It comprises the points $x\in X$
that are not smooth, and puts a scheme structure on this locus, which is usually
nonreduced.

\begin{proposition}
\mylabel{jacobian subscheme}
 Suppose that the nonsmooth locus of $X$ is zero-dimensional,
that is, the jacobian subscheme $X'\subset X$ is zero-dimensional.
Then the number $h^0(X,\O_{X'})$ is a multiple of $p$.
\end{proposition}

\proof
Follows as above from Lemma \ref{multiple p}.
\qed

\medskip
If $x\in X$ is an isolated nonsmooth point, then the local length $l(\O_{X',x})$
 is   called the \emph{Tjurina number}
of the singularity $x\in X$.
We thus see that the $p$-divisibility of   Tjurina numbers  
is a necessary condition for an $F$-scheme  to descend to 
regular schemes.

\section{Finite projective dimension}
\mylabel{finite projective dimension}

\medskip
We keep the assumptions from the preceding section, such that
$X$ is an $F$-scheme of finite type that descends to a regular $E$-scheme $Y$,
where $E\subset F$ is a purely inseparable field extension.
The goal of this section is to exploit finer properties of coherent sheaves $\shF_X=\shF_Y\otimes_E F$
related to commutative algebra rather than field theory.

Suppose $\shM$ is a coherent $\O_X$-module.
Given a point $x\in X$,
the \emph{projective dimension} $\pd(\shM_x)$ is defined as the infimum over all numbers $n$
such that there exists a finite free resolution
$$
0\lra F_n\lra F_{n-1}\lra\ldots\lra F_1\lra F_0\lra \shM_x\lra 0,
$$
where $F_n$ are free $\O_{X,x}$-modules of finite rank.
In case that no such resolution exists, one writes $\pd(\shM_x)=\infty$.

\begin{lemma}
\mylabel{projective dimension}
Let $\shF_Y$ be a coherent sheaf on $Y$, and $\shF_X=\shF_Y\otimes_E F$
the induced sheaf on $X$. Then $\shF_{X,x}$ has finite projective dimension
for all $x\in X$.
\end{lemma}

\proof
Let $y\in Y$ be the image of $x\in X$. Since $\O_{Y,y}$ is regular by assumption,
the stalk $\shF_{Y,y}$ has finite projective dimension, according to Serre's Criterion.
Since $X\ra Y$ is flat, this implies that $\shF_{X,x}$ has finite projective dimension as well.
\qed

\medskip
In particular, the stalks of $\Omega^1_{X/F}$ have finite projective dimension.
From this, we immediately deduce the following   facts:

\begin{proposition}\mylabel{kaehler differentials}
Any coherent $\O_X$-module obtained from $\Omega^1_{X/F}$ by 
taking tensor powers, symmetric powers, alternating powers,  
Fitting ideals, duals, biduals, kernels and cokernels of biduality maps etc.\
have stalks with finite projective dimension.
\end{proposition}

In particular, the stalks of the jacobian ideal $\shJ\subset\O_X$ have finite
projective dimension. 
In case the jacobian subscheme $X'\subset X$ is zero-dimensional, we have
as an invariants the  Tjurina numbers, that is, the lengths $l(\O_{X,x}/\shJ_x)\geq 0$ of the jacobian subscheme at $x\in X$.

Recall that for any coherent ideal $\shI\subset\O_X$, we define the \emph{bracket ideal}
$\shI^{[p]}\subset\O_X$ as the ideal whose stalks are generated by  $f^p$, $f\in\shI_x$.

\begin{theorem}
\mylabel{length criterion}
Suppose that the jacobian subscheme $X'\subset X$ is zero-dimensional, and let $d=\dim(X)$.
Then for all closed points $x\in X$, the length formula
$l(\O_{X,x}/\shJ_x^{[p]})=p^d l(\O_{X,x}/\shJ_x)$ holds.
\end{theorem}

\proof
The scheme $X$ is locally of complete intersection by Proposition \ref{complete intersection}.
In case $x\in X$ is smooth, both lengths in question are zero, so it suffices to
treat only the case $x\in X'$.
Since $X'$ is discrete, the ideal $\shJ_x\subset\O_{X,x}$ is
$\maxid_x$-primary.
According to Miller's result \cite{Miller 2003}, Corollary 5.2.3,  the length formula 
$l(\O_{X,x}/\shI_x^{[p]})=p^d l(\O_{X,x}/\shI_x)$ holds for an ideal $\shI_x\subset\O_{X,x}$
defining a zero-dimensional subscheme
if and only if  $\shI_x$ has finite projective dimension.
Whence Proposition \ref{kaehler differentials} implies the assertion.
\qed

\medskip

\begin{example}
\mylabel{nondescent}
Assume that $X$ is a normal surface over   in characteristic $p=2$.
Suppose $x\in X$ is a rational point so that the local ring $\O_{X,x}$ 
is a rational double point of type $E_8^3$. If $k\subset\bar{k}$ is an algebraic closure,
we have $\O_{X,\bar{x}}^\wedge\simeq k[[x,y,z]]/(f)$, were $f=z^2+x^3+y^5+y^3z$,
as explained in Artin's paper \cite{Artin 1977}.
We have
$$
\shJ_{\bar{x}}=(x^2,y^4+y^2z,y^3)=(x^2,y^2z,y^3).
$$
From this one easily computes the lengths  $l(\O_{X,x}/\shJ_x)=10$ and
$l(\O_{X,x}/\shJ_x^{[2]})=44$, either by hand or with a computer. 
This implies that the $F$-scheme $X$ does not descend to an $E$-scheme that is regular.
Note that the local fundamental group of the strict henselization of $x\in X$ is
trivial, so Proposition  \ref{strictly henselian} does not allow this conclusion.
\end{example}

Let us finally have a closer look the tangent sheaf $\Theta_X=\shHom(\Omega^1_{X/k},\O_X)$.
Clearly, $\Theta_X=\Theta_Y\otimes_E F$ holds, so the stalks of $\Theta_X$ have finite projective
dimension, according to Proposition \ref{kaehler differentials}.
One can say more in small codimensions:

\begin{proposition}
\mylabel{tangent sheaf}
Let $x\in X$ be a point of codimension two. 
Then the stalk $\Theta_{X,x}$ is a free $\O_{X,x}$-module.
\end{proposition}

\proof
The scheme $X$ is Cohen--Macaulay, whence $\depth(\O_{X,x})=2$.
Being a dual,  the sheaf $\Theta_X$ satisfies Serre's Condition $(S_2)$, according
to \cite{Hartshorne 1994}, Theorem 1.9. In particular, $\depth(\Theta_{X,x})\geq 2$.
We stalks of $\Theta_X$ have finite projective dimension, whence the Auslander--Buchsbaum Formula
$$
\pd(\Theta_{X,x})+\depth(\Theta_{X,x})=\depth(\O_{X,x}),
$$
holds. We infer   $\pd(\Theta_{X,x})=0$. This means that the stalk $\Theta_{X,x}$
is free.
\qed

\begin{remark}
 \mylabel{irreducible codimension}
Lipman (\cite{Lipman 1965}, Proposition 5.2)
proved that if $\Theta_{X}$ is locally free, than all irreducible components
of the jacobian subscheme $X'\subset X$ have codimension $\leq 2$.
\end{remark}

\begin{remark}
\mylabel{length free}
 Let me point out that there is a close connection between the length formula
of Theorem \ref{length criterion} and the local freeness of the tangent sheaf
in Proposition \ref{tangent sheaf}, which is independent of our assumption
that $X$ descends to a regular scheme.
Suppose that $S$ is the spectrum of $F[x,y,z]/(f)$ for some polynomial $f$,
and assume that $S$ is normal.
As explained in \cite{Artin 1976}, Theorem 6.2, dualizing the
short exact sequence $0\ra\shI/\shI^2\ra\Omega^1_{\AA^3_E/E}\mid_S\ra\Omega^1_S\ra 0$
yields an exact sequence
$$
0\lra\Theta_S\lra\Theta_{\AA^3_F}|_S\lra(\shI/\shI^2)^\vee\lra\shT_S^1\lra 0,
$$
where $\shI=\O_{\AA^3_F}(-S)$, and the
 cokernel on the right $\shT^1_S$ is isomorphic to the
sheaf $\shExt^1(\Omega^1_{S/F},\O_S)$ of first-order extension.
The preceding exact sequence shows that  the annulator ideal of $\shT^1_S$ 
coincides with the jacobian ideal $\shJ$, and that $\shT^1_S$ is an invertible as
 module over $\O_S/\shJ$. It follows that
$\pd(\Theta_S)=\pd(\shJ)-2$,
and that $\Theta_S$ is locally free if and only if the length formula
$l(\O_{S,s}/\shJ_s^{[p]})=p^2 l(\O_{S,s}/\shJ_s)$
holds. Note that the latter formula is fairly easy to check, either by hand or
with a computer. In contrast, finding an explicit basis of the stalks of 
the tangent sheaf seems to involve some nontrivial guess-work.
\end{remark}

\section{Invertible summands of   cotangent sheaf}
\mylabel{invertible summands cotangent}

Let $X$ be an $F$-scheme of finite type, and let  $E\subset F$  be a subfield so that the field 
extension $E\subset F$ is   purely inseparable. In the preceding sections, we saw necessary conditions for $X$
to descend to a regular $E$-scheme,
involving local fundamental groups,    degrees of residue fields,   jacobian ideals,
and the tangent sheaf. These conditions, however, are   not always strong enough do
rule out that $X$ descends to a regular scheme.
In this section we give another conditions in terms of the sheaf of K\"ahler differentials
$\Omega^1_{X/F}$.
Together with the other criteria, this 
will be enough to settle the case of rational double points, which   will
happen in the last section.

Let me start with the following simple observation from commutative algebra:
Suppose  $R$ is a local    ring.
Let $M$ be an $R$-module, $M^\vee=\Hom(M,R)$ the dual module,
and $\Phi: M\times M^\vee\ra R$ the canonical bilinear map $(a,f)\mapsto f(a)$.

\begin{lemma}
\mylabel{evaluation pairing}
The following conditions are equivalent:
\renewcommand{\labelenumi}{(\roman{enumi})}
\begin{enumerate}
\item
The induced linear map $\Phi: M\otimes M^\vee\ra R$ is surjective.
\item There are elements $a\in M$ and $f\in M^\vee$ with $f(a)\not\in\maxid_R$.
\item There is an $R$-module $M'$ and a decomposition $M\simeq R\oplus M'$.
\end{enumerate}
\end{lemma}

\proof
The implications (i)$\Rightarrow$(ii) and (iii)$\Rightarrow$(i) are easy.
The task is to verify (ii)$\Rightarrow$(iii).
Suppose   $f(a)\not\in\maxid_R$ for some $a\in M$ and $f\in M^\vee$.
Then $f(a)\in R$ must be a unit, because $R$ is local.
It follows that the map $M\ra R$, $x\mapsto f(x)$ is surjective.
Whence there is an exact sequence
$$
0\lra M'\lra M\stackrel{f}{\lra} R\lra 0
$$
with $M'=\ker(f)$. Such a sequence splits because $R$ is free.
\qed

\medskip
We say that an $R$-module $M$ \emph{has an invertible   summand}  if the
equivalent conditions of the lemma hold.
In the noetherian situation, we may check this by passing back and forth  to formal completions.
More generally:

\begin{proposition}
\mylabel{formal completion}
Suppose $M$ is of finite presentation, and let $R\ra R'$ be a faithfully flat ring homomorphism.
Then the $R$-module $M$ has an invertible   summand if and only this holds for the induced $R'$-module $M'=M\otimes_R R'$
\end{proposition}

\proof
The condition is clearly necessary. Conversely, suppose
that the evaluation map $M'\otimes(M')^\vee\ra R'$ is surjective.
Using that the canonical homomorphism $M^\vee\otimes_R R'\ra(M\otimes_R R')^\vee$ is bijective, 
we infer that $M\otimes M^\vee\ra R$ is surjective as well.
\qed

\medskip
We shall apply this concept to the stalks $M=\Omega^1_{X,x}$ of the cotangent sheaf:

\begin{theorem}
\mylabel{direct summands}
Suppose our $F$-scheme of finite type $X$ is equidimensional, of dimension $\dim(X)\geq 1$, and  descends to a regular scheme.
Then for all points $x\in X$ of codimension $\leq 2$, the stalk $\Omega^1_{X,x}$
has an  invertible   summand.
\end{theorem}

\proof
If $x\in X$ is smooth, then $\Omega^1_{X/k}$ is locally free at $x\in X$ of rank $\dim(X)\geq 1$,
hence there is nothing to show.
Assume that $x\in X$ is not smooth. Let $y\in Y$ be its image.
It suffices to check that $\Omega^1_{Y/E,y}$ has an  invertible summand.
Whence our task   is to find a local vector field $\delta\in\Theta_{Y,y}$ together
with a local section $s\in\O_{Y,y}$ with $\delta(s)=1$.

First of all, we may assume that $Y$ is affine.
Let $Y'=\overline{\left\{y\right\}}$ be the reduced closure of the point $y\in Y$,
and $L=\kappa(y)$ be its residue field. Then the finitely generated field extension
$E\subset L$ has transcendence degree $\trdeg_E(L)=\dim(Y')$.
Let $E\subset K\subset L$ be the separable closure of $E(x_1,\ldots,x_m)$,
where $x_1,\ldots,x_m\in L$ is a transcendence basis.
Then the finite field extension $K\subset L$ is purely inseparable.
According to Proposition \ref{linearly disjoint}, we have $K\neq L$.
It follows that $\Omega^1_{L/K}\neq 0$, by \cite{A 4-7}, \S6, No.\ 3, Theorem 3.
In light of the exact sequence
$\Omega^1_{L/E}\ra \Omega^1_{L/K}\ra 0$,
there must be  some $\lambda\in L$ with nonzero differential $d\lambda\in\Omega^1_{L/E}$.
After replacing $Y$ by a suitable affine open subscheme, we may extend
 the value $\lambda\in L=\kappa(y)$ to some  section $s'\in H^0(Y,\O_{Y'})$.
By the same token, we may furthermore assume that the coherent $\O_{Y'}$-module
$\Omega^1_{Y'/E}$ is locally free.
Moreover, we may lift $s'$ to a section $s\in H^0(Y,\O_Y)$. Then $ds(y)\neq 0$
as element of the $\kappa(y)$-vector space $\Omega^1_{Y/E}(y)$.

Consider the biduality map $\Omega^1_{Y/E}\ra(\Omega^1_{Y/E})^{\vee\vee}$, which sends
differentials $\sum s_idf_i$ to the evaluation map $\delta\mapsto\delta(\sum s_idf_i)$.
Let us denote this evaluation map by the symbol $(\sum s_i df_i)^{\vee\vee}$,
which is thus a local section of $\Theta_Y^\vee$.
Our aim now is to verify that the value at $y$ of the evaluation map $(ds)^{\vee\vee}$, which is 
an element of  $\Theta_{Y/E}^\vee(y)$, does not vanish. 
Here the problem is that the biduality map for K\"ahler differentials on $Y$
is not necessarily bijective. However, we know that the biduality map for K\"ahler differentials
on $Y'$ is bijective. Before we exploit this, let me make a brief digression on biduality maps:

Suppose $\shF$ and $\shG$ are coherent sheaves on $Y$ and $Y'$, respectively.
Set $\shF'=\shF\otimes\O_{Y'}$, and assume we have a   homomorphism $f:\shF'\ra\shG$.
Consider the canonical restriction map
$\shHom(\shF,\O_Y)\otimes\O_{Y'}\ra\shHom(\shF',\O_{Y'})$.
Applying it twice, we get a canonical map
$\shHom(\shHom(\shF,\O_Y),\O_Y)\otimes\O_{Y'}\ra \shHom(\shHom(\shF',\O_{Y'}),\O_{Y'})$.
We may compose this with the map 
induced from $f:\shF'\ra\shG$ to obtain a map
$$
\varphi:\shHom(\shHom(\shF,\O_Y),\O_Y)\lra \shHom(\shHom(\shG,\O_{Y'}),\O_{Y'}).
$$
This yields  a commutative diagram
$$
\begin{CD}
\shF @>>> \shHom(\shHom(\shF,\O_Y),\O_Y)\\
@VfVV @VV\varphi V\\
\shG @ >>> \shHom(\shHom(\shG,\O_{Y'}),\O_{Y'}),
\end{CD}
$$
where the horizontal maps are the biduality maps on $Y$ and $Y'$, respectively.

We may apply this to the sheaves $\shF=\Omega^1_{Y/E}$ and $\shG=\Omega^1_{Y'/E}$
and the canonical surjection $\Omega^1_{Y/E}\otimes\O_{Y'}\ra\Omega^1_{Y'/E}$.
This gives us a commutative diagram
$$
\begin{CD}
\Omega^1_{Y/E} @>>> \shHom(\shHom(\Omega^1_{Y/E},\O_Y),\O_Y)\\
@VVV @VVV\\
\Omega^1_{Y'/E} @ >>> \shHom(\shHom(\Omega^1_{Y'/E},\O_{Y'}),\O_{Y'}).
\end{CD}
$$
The lower horizontal map is bijective, because $\Omega^1_{Y'/E}$ is locally free on $Y'$.
Whence $(ds')^{\vee\vee}$ does not vanish at $y\in Y$.
It follows that $(ds)^{\vee\vee}$ does not vanish at $y$ as well.

Now recall that the stalk   $\Theta_{Y,y}$ is free,
by Proposition \ref{tangent sheaf}. Then $(\Omega^1_{Y/E})^{\vee\vee}_y$ is free as well.
Set  $e_1=(ds)^{\vee\vee}$ and extend it
to a basis   $e_1,\ldots,e_n\in(\Omega^1_{Y/E})^{\vee\vee}_y$.
 Let $\delta_1,\ldots,\delta_n\in\Theta_{Y,y}$  be the corresponding  
dual basis.  Then $\delta=\delta_1$ is the desired local vector field:
We have $\delta(s)=\delta(ds)=\delta(ds)^{\vee\vee}=1$
by construction.
\qed

\medskip
It follows from Lipman's results (\cite{Lipman 1965}, Proposition 8.1)   that
$\Omega^1_{X/F}$ is torsion free if and only if our scheme $X$ is geometrically normal.
According to \cite{Evans; Griffith 1985}, Theorem 2.14
any torsion free module of finite type over a local noetherian ring
is an extension of an ideal by a free module. Such extensions, however, do not necessarily split.
So the preceding result puts a nontrivial condition on the stalk $\Omega^1_{X/F,x}$.

It turns out that the result is well-suited to rule out that certain $F$-schemes 
descend to regular schemes.
Suppose $X$ is 2-dimensional, normal, and locally of complete intersection.
Let $x\in X$ be a closed point, with embedding dimension $\edim(\O_{X,x})=3$.
Write $\O_{X,x}^\wedge=\Spec F[[u,v,w]]/(f)$ for some nonzero polynomial $f$.
\comment{power series instead?}
Using   notation from physics, we write $f_u=\frac{\partial f}{\partial u}$ etc.\ for   partial derivatives.

\begin{corollary}
\mylabel{partial derivatives}
Assumptions as above. Suppose also that $X$ descends to a regular  scheme.
Then, after a suitable permutation of the indeterminates $u,v,w$, the ideal $\ida=(f_u,f_v,f)$ inside
the formal power series ring $F[[u,v,w]]$ is a parameter ideal, and furthermore $f_w\in\ida$ holds.
\end{corollary}

\proof
We may replace   $X$ by  
  $\Spec(F[u,v,w]/(f))\subset\AA^3_F$, such that our point $x\in X$ corresponds to the origin $0\in\AA^3_F$.
Now we have an exact sequence
$$
0\lra\shI/\shI^2\lra\Omega^1_{\AA^3_F/F}|_X\lra\Omega^1_{X/F}\lra 0,
$$
where $\shI=\O_{\AA^3_F}(-X)$.
Dualizing, we get an exact sequence
$$
(\shI/\shI^2)^\vee\longleftarrow \Theta_{\AA^3_F}|_X\longleftarrow\Theta_X\longleftarrow 0.
$$
Whence any any local vector field $\delta\in\Theta_{X,x}$ is the restriction of 
some local vector field of the form 
$$
\tilde{\delta}=rD_u+sD_v+tD_w,
$$
subject to the condition $\tilde{\delta}(f)\in(f)_\maxid$, where $\maxid=(u,v,w)$.
Here $D_u$ etc.\ denotes the derivation given by taking partial derivative.
The coefficients are local sections $r,s,t\in\O_{\AA^3,x}=k[u,v,w]_\maxid\subset F[[u,v,w]]$.

By Theorem \ref{direct summands}, the stalk $\Omega^1_{X,x}$ contains an invertible direct summand.
Whence there is a derivation $\tilde{\delta}=rD_u+sD_v+tD_w$ so that first 
$\tilde{\delta}(f)=rf_u+sf_v+tf_w\in (f)_\maxid$, and
second that at least on of the local sections
$\tilde{\delta}(u)=r$ and $\tilde{\delta}(v)=s$ and $\tilde{\delta}(w)=t$ is invertible inside   $ F[u,v,w]_\maxid$. After a permutation of indeterminates, we may
assume that $\tilde{\delta}(w)=t$ is invertible.
From this we infer that $f_w\in (f_u,f_v,f)$  
holds, where the ideals are considered inside the local ring $k[u,v,w]_\maxid$, or equivalently in $F[[u,v,w]]$.

It remains to check that $(f_u,f_v,f)$ is a parameter ideal.
Indeed: The ideal $(f_u,f_v,f_w,f)$ defines the nonsmooth locus of $X$.
This locus is discrete, because $X$ is a normal surface by assumption.
It follows that $(f_v,f_w,f)=(f_u,f_v,f_w,f)$ is a parameter ideal in $F[[u,v,w]]$.
\qed

\section{Rational double points}
\mylabel{rational double points}

The goal of this section is to determine which rational double points
descend to regular schemes.
Suppose that $F$ is an algebraically closed ground field of characteristic $p>0$, and let $X$ be a normal
surface over $F$. A singular point $x\in X$ is called a \emph{rational double point}
if the singularity is rational and Gorenstein (equivalently: rational and of multiplicity two). 
Throughout this section, we tacitely assume that the field $F$ is not algebraic over its prime field.

Rational double points are classified according to the intersection graph of
the exceptional divisor on the minimal resolution of singularities; these intersection graphs
correspond to the simply laced Dynkin diagrams $A_n$, $n\geq 1$ and $D_n$, $n\geq 4$
and $E_6$, $E_7$, $E_8$. In characteristic zero, two rational double points are formally isomorphic
if and only if they have the same Dynkin type. According to Artin \cite{Artin 1977},
this is not true   in characteristic $p=2,3,5$.
However, there are still only finitely many formal isomorphism classes, and Artin gave a list
of formal equations for these isomorphism classes. Before we come to this,
let my reprove the following result, which is due to Hirokado \cite{Hirokado 2004}:

\begin{theorem}
\renewcommand{\labelenumi}{(\roman{enumi})}
\begin{enumerate}
\item
A rational double point of type $A_n$, $n\geq 1$ descends to a regular scheme if and only if $n+1=p^{e}$ for some exponent $e\geq 1$.
\item
In characteristic $p\geq 3$,   rational double points of type $D_n$   do not
descend to regular schemes
\item For  $p\geq 7$, rational double points of type $E_n$ do not
descend to regular schemes.
\end{enumerate}
\end{theorem}

\proof
Rational double points of type $A_n$ are defined by   $f=z^{n+1}-xy$.
If we have $n+1=p^e$, then this rational double point descends to regular schemes according
to Proposition \ref{regular singular}. Conversely, suppose that $n+1$ is not a $p$-power. Then the local Picard group $\Pic^\loc$, which is cyclic of order $n+1$,
is not a $p$-group. By Corollary \ref{class group}, such a rational double point does not appear on 
generic geometric fibers. 
The local Picard group of a rational double point of type $D_4$ is of order four.
Using Corollary \ref{class group} again, we conclude that for $p\geq 3$ rational double points of type $D_n$ do not
appear on geometric generic fibers. This settles assertion (i) and (ii).
$$
\begin{array}[t]{l|c|c|c|c|c}
\text{rational double point}           &  A_n  &  D_n   & E_6 & E_7 & E_8 \\ 
\hline&&&&&\\[-2ex] 
\text{$\pi_1^\loc$ for $p=0$} & \text{cyclic}   & \text{dihedral} &  \tilde{A}_4 & \tilde{S}_4  &  \tilde{A}_5 \\
\hline&&&&&\\[-2ex]  
\text{order of $\pi_1^\loc$ for $p=0$} & n+1   & 2(n-2) &  24 & 48  &  120 \\
\hline&&&&&\\[-2ex]
\text{order of $\Pic^\loc$} & n+1 & 4 & 3 & 2 & 1\\
\end{array}
$$

\vspace{1em}
\centerline{Table \stepcounter{figure}\arabic{figure}: Local fundamental and Picard groups of rational double points}
\vspace{1em}

It remains to treat the $E_n$ case.
According to \cite{Artin 1977}, Proposition 2.7, the  local fundamental group $\pi_1^\loc$ of  rational double points
in characteristic $p\geq 7$ 
are  \emph{tame}. This simply means that they have  order prime to the characteristic.
It follows that $\pi_1^\loc$ is the maximal prime-to-$p$ quotient of the corresponding local fundamental group
in characteristic zero. Their orders appear in Table 1 below (confer, for example, \cite{Lamotke 1986}).
We infer that $\pi_1^\loc\neq 0$ in characteristic  $p\geq 7$.
Proposition \ref{strictly henselian} tells us that the rational double points of type $E_n$ for $p\geq 7$ do not
appear on geometric generic fibers.
\qed

\medskip
Let us now turn to rational double points of type $D_n$, $n\geq 4$ in characteristic two.
According to \cite{Artin 1977}, they are
subdivided into $\lfloor n/2\rfloor$ isomorphism classes $D_n^r$, depending on an additional integral parameter $0\leq r \leq \lfloor n/2\rfloor -1$. The formal equations
$f(x,y,z)=0$ are as follows:
\begin{equation}
\label{dn singularities}
\begin{split}
D_{2m}^0:\quad & f= z^2+x^2y +xy^m, \\
D_{2m}^r:\quad &f=z^2+x^2y +xy^m + xy^{m-r}z,\\
D_{2m+1}^0:\quad & f=z^2+x^2y +y^mz,\\
D_{2m+1}^r:\quad & f=z^2+x^2y +y^mz +xy^{m-r}z.
\end{split}
\end{equation}

\begin{remark}
Note that the polynomial $f=z^2+x^2y +xy^m$ defining $D_{2m}^0$ is contact equivalent to 
$g=z^2+x^2y +xy^m + xy^mz$, which would be the case $r=0$ in $D_{2m}^r$. 
This fact follows easily from the methods in \cite{Roczen 1996}, Section 2. Recall that 
\emph{contact equivalence} means that the two equations define formally isomorphic singularities.
The sitution for $D_n$ with $n$ odd is similar. We thus see that it is not really necessary to distinguish the cases $r=0$ and $r>0$ when it comes to giving the formal equations.
\end{remark}

\begin{theorem}
\mylabel{d4 p=2}
In characteristic $p=2$, a rational double point of type $D_n^r$
descends to a regular scheme if and only if  $r=0$.
\end{theorem}

\proof
The condition is sufficient: Rational double points of type $D_{2m}^0$, which are formally given
by $f= z^2+x^2y +xy^m$, descend to regular schemes according to Proposition \ref{regular singular}.

The case $D_{2m+1}^0$ is slightly more complicated, due to the appearance of the   linear $z$-monomial in the defining polynomial  
$f=z^2+x^2y +y^mz$. Choose a nonperfect subfield $E\subset F$ and an element $\lambda\in E$
that is not a square. Consider the spectrum of $A=E[x,y,z]/(f+\lambda)$. The ring $A$ is a flat $E[x,y]$-algebra
of degree two. The fiber ring $A/(x,y)A\simeq E(\sqrt{\lambda})$ is regular, which implies  that the ring $A$ is regular.
It remains to check that $f$ and $f+\lambda$ define isomorphic singularities over the algebraically closed field $F$.
To do so, we construct an automorphism of the ring $ F[[x,y]][z]$ sending $f+\lambda$ to $f$.
Clearly, the substitution $z\mapsto z+\sqrt{\lambda}$ sends $f+\lambda$ to $f+\sqrt{\lambda}y^m$.
If $m=2l+1$ is odd, then the additional substitution $x\mapsto x+\lambda^{1/4}y^l$ achieves our goal.
If the number $m=2l$ is even, we make the substitution $z\mapsto z+\lambda^{1/4}y^l$
instead, which maps $f+\sqrt{\lambda}y^m$ to $f+\lambda^{1/4}y^{m+m/2}$. Repeating substitutions of the latter kind,
we are, after finitely many steps, in the position to apply a substitution of the former kind.

The condition is also necessary: Set $S=\Spec F[x,y,z]/(f)$, where $f$ is the equation of a rational double point
of type $D_n^r$ with $1\leq r\leq\lfloor n/2\rfloor -1$, as given in (\ref{dn singularities}). Let $s\in S$ be the singular point.
We shall verify that $\Omega^1_{S,s}$ contains no invertible summand. To do so, we apply Corollary 
\ref{partial derivatives}, which
reduced our problem at hand to a calculation involving the partial derivatives $f_x,f_y,f_z$.
Suppose first that $n=2m$ is even.
We compute
\begin{equation*}
\begin{split}
f&=z^2+x^2y +xy^m + xy^{m-r}z, \\
f_x&=y^m+y^{m-r}z,\\
f_y&=x^2+mxy^{m-1}+(m-r)xy^{m-r-1}z,\\
f_z&=xy^{m-r}.
\end{split}
\end{equation*}
Consider the ideal  
$
I_x=(f_y,f_z,f)$ 
inside the polynomial ring $F[x,y,z]$. Reducing modulo $x$  we have $I_x\cong(z^2)$, and  see that the induced ideal  
 $I_x\O_{S,s}$ has height one.
An analogous argument applies to the ideal $I_y=(f_x,f_z,f)$, where we reduce modulo $y$.
Finally, consider the ideal $I_z=(f_x,f_y,f)$. Computing modulo $z$ we have
$$
I_z=(f_x,f_y,f)\cong (y^m,x^2+mxy^{m-1},x^2y)\cong(y^m,x^2+mxy^{m-1}).
$$
It follows that the residue classes of $x^iy^j$ with $0\leq i\leq 1$ and $0\leq j\leq m-1$ form an
$F$-vector space basis of the residue class ring $R=k[x,y]/(y^m,x^2+mxy^{m-1})$.
Obviously, the residue class of $f_z=xy^{m-r}$ is nonzero, and the Artin ring $R$ is local.
This implies $f_z\not\in I_z\O_{S,s}$.
Invoking Corollary \ref{partial derivatives}, we deduce that the rational double points of type $D_{2m}^r$, $r>0$ do not
appear on geometric generic fibers.

The case that $n=2m+1$ is odd can be treated with similar arguments. Here we have 
\begin{equation*}
\begin{split}
f&=z^2+x^2y +y^mz + xy^{m-r}z, \\
f_x&=y^{m-r}z,\\
f_y&=x^2+my^{m-1}z+(m-r)xy^{m-r-1}z,\\
f_z&=y^m+xy^{m-r}.
\end{split}
\end{equation*}
We proceed as above: Computing modulo $x$, we see
$$
I_x=(f_y,f_z,f)\cong( my^{m-1}z,y^m,z^2),
$$
and   infer that the residue class of the partial derivative $f_x=y^{m-r}z$ inside the local ring 
$R=k[y,z]/(my^{m-1}z,y^m,z^2 )$
is nonzero, provided $r\geq 2$ or $m$ even. This implies $f_x\not\in I_x\O_{S,s}$.
It remains to check the case $r=1$ and $m$ odd. Then  we compute modulo $y^m,xy$ and have
$I_x\cong(x^2+y^{m-1}z,z^2)$, and easily check that $f_x=y^{m-1}z$ is nonzero inside
$R=k[x,y,z]/(x^2+y^{m-1},z^2,y^m,xy)$. Again $f_x\not\in I_x\O_{S,s}$.

Computing modulo $y$, we see
$
I_y=(f_x,f_z,f)\cong(z^2),
$
such that $I_y\O_{S,s}$ has height one.
Finally, we have $I_z\cong 0$  modulo $x,z$, such that  $ I_z\O_{S,s}$ has height one.
According to Corollary \ref{partial derivatives}, the rational double points of type $D_{2m+1}^r$ with $r\geq 1$ do not
appear on geometric generic fibers.
\qed

\medskip
We finally come to the rational double points of type $E_6,E_7,E_8$, which are the most challenging cases:

\begin{theorem}
\renewcommand{\labelenumi}{(\roman{enumi})}
\begin{enumerate}
\item
Suppose   $p=5$. Among the four rational double points of type $E_n$, only $E_8^0$ descends to a regular scheme.
\item
Suppose $p=3$. Among the seven rational double points of type $E_n$, only $E_6^0$ and $E_8^0$ descend to regular schemes.
\item
Suppose $p=2$. Among the eleven rational double points of type $E_n$, only $E_7^0$ and
$E_8^0$ descend to regular schemes.
\end{enumerate}
\end{theorem}

\proof
The formal equations $f(x,y,z)=0$ for the rational double points were determined by Artin \cite{Artin 1977},
and appear in Tables 2--4.
Proposition \ref{regular singular} 
immediately tells us  that the rational double points of type $E_n^0$ as stated in the assertion
descend to regular schemes.
The task is to argue that the remaining rational double points do not.

To do so, I have collected further information   in the tables.
The third columns give the local fundamental groups, which were determined by Artin \cite{Artin 1977}.
The forth columns contain
informations on the jacobian ideal $J\subset F[x,y,z]/(f)$ and the corresponding bracket ideal $J^{[p]}$.
Recall that $J$, $J^{[p]}$ are induced from  the ideals $(f_x,f_y,f_z,f)$, $(f_x^p,f_y^p,f_z^p,f)$ inside the polynomial ring
$F[x,y,z]$, respectively. Column three contains the length of the quotient by $J$ and $J^{[p]}$.
According to Theorem \ref{length criterion}, the length formula 
$$
l(F[x,y,z]/(f_x,f_y,f_z,f)) = p^2l(F[x,y,z]/(f_x^p,f_y^p,f_z^p,f))
$$
holds if the singularity descends to a regular surface. By Remark \ref{length free},
the length formula holds if and only if
the tangent sheaf is locally free. The latter information appears in the last column of the tables. 
Using  the information on the local Picard group for Proposition \ref{local fundamental}, and the information about the length for
Theorem \ref{length criterion}, we rule out all   rational double points, except the following three in characteristic $p=2$:
\begin{equation}
\label{en singularities}
\begin{split}
E_7^1:\quad & f= z^2+x^3+xy^3+x^2yz, \\
E_7^2:\quad & f=z^2+x^3+xy^3+y^3z,\\
E_8^1:\quad & f=z^2+x^3+y^5+xy^3z.
\end{split}
\end{equation}
For them, the local fundamental group is trivial and the tangent sheaf is locally free.
We shall descard them   by showing that the cotangent sheaf does not contain an invertible direct
summand at the singularity. As in the case of $D_n$-singularities, we will apply Corollary \ref{partial derivatives}.
Let me treat the case of $E_7^1$. The partial derivatives are
\begin{equation*}
\begin{split}
f&=f= z^2+x^3+xy^3+x^2yz, \\
f_x&=x^2+y^3,\\
f_y&=xy^2+x^2z,\\
f_z&=x^2y.
\end{split}
\end{equation*}
Consider the ideal $I_x=(f_y,f_z,f)$ inside the power series ring $F[[x,y,z]]$.
We have $I_x\subset (x,z)$, whence $I_x$ is not a parameter ideal.
Next, consider the ideal 
$$
I_y=(f_x,f_z,f)=(x^2+y^3,x^2y,z^2)=(x^2+y^3,y^4,z^2).
$$
The residue classes of the monomials $x^iy^jz^k$ with $0\leq i,k\leq 1$ and $0\leq j\leq 3$
constitute an $F$-vectors space basis for the quotient ring $F[[x,y,z]]/I_y$.
It follows that the residue class of $f_y=xy^2+x^2z$ is nonzero, whence $f_y\not\in I_y$.
Finally, consider the ideal $I_z=(f_x,f_y,f)=(x^2+y^3,xy^2+y^3z,z^2+y^4z)$.
Let us compute modulo $z$: Then $I_z\cong (x^2+y^3,xy^2)=(x^2+y^3,xy^2,y^5)$,
and it is straightforward to see that the residue classes of $1,y,\ldots,y^4,x,xy$ form
a $F$-vector space basis modulo $I_z+(z)$. Therefore, $f_z=x^2y\cong y^4$ is not contained in $I_z$.
Using Corollary \ref{partial derivatives}, we conclude that rational double points of type $E_7^1$ do not
appear on geometric generic fibers.
The remaining two cases $E_7^2$ and $E_8^1$ are somewhat similar, and left to the reader.
\qed

$$
\begin{array}[t]{l|l|l|l|l}
\text{}        &  \text{formal relation}   &  \pi_1^\loc &   l(\O/\shJ),l(\O/\shJ^{[2]}) & \text{$\Theta$ free} \\
\hline&&&&\\[-2ex]
E_6^0 & z^2+x^3+y^2z       & C_3 & 8,32 & \text{yes}  \\
E_6^1 & z^2+x^3+y^2z+xyz   & C_6 & 6,28 & \text{no} \\
\hline&&&&\\[-2ex]
E_7^0 & z^2+x^3+xy^3       &0    & 14,56 & \text{yes}  \\
E_7^1 & z^2+x^3+xy^3+x^2yz &0 & 12,48 & \text{yes} \\
E_7^2 & z^2+x^3+xy^3+y^3z  & 0 & 10,40 & \text{yes} \\
E_7^3 & z^2+x^3+xy^3+xyz   & C_4 & 8,35 & \text{no} \\
\hline&&&&\\[-2ex]
E_8^0 & z^2+x^3+y^5        & 0 & 16,64 & \text{yes} \\
E_8^1 & z^2+x^3+y^5+xy^3z  & 0 & 14,56 & \text{yes} \\
E_8^2 & z^2+x^3+y^5+xy^2z  & C_2 & 12, 48 & \text{yes} \\
E_8^3 & z^2+x^3+y^5+y^3z   & 0 & 10,44 & \text{no} \\
E_8^4 & z^2+x^3+y^5+xyz    & C_3\rtimes C_4 & 8,37 & \text{no} 
\end{array}
$$

\vspace{1em}
\centerline{Table \stepcounter{figure}\arabic{figure}: Rational double points of  type $E_n$ in characteristic $p=2$.}
\vspace{1em}

$$
\begin{array}[t]{l|l|l|l|l}
\text{}         &  \text{formal relation}   &  \pi_1^\loc &   l(\O/\shJ),l(\O/\shJ^{[3]}) & \text{$\Theta$ free}\\
\hline&&&&\\[-2ex]
E_6^0 & z^2+x^3+y^4        & 0 & 9,81 &  \text{yes} \\
E_6^1 & z^2+x^3+y^4+x^2y^2 & C_3 & 7,71 &  \text{no}\\
\hline&&&&\\[-2ex]
E_7^0 & z^2+x^3+xy^3       & C_2   & 9,81 &  \text{yes} \\
E_7^1 & z^2+x^3+xy^3+x^2y^2 &C_6 & 7,75 &   \text{no}\\
\hline&&&&\\[-2ex]
E_8^0 & z^2+x^3+y^5        & 0 & 12,108 &   \text{yes}\\
E_8^1 & z^2+x^3+y^5+x^2y^3  & 0 & 10,99 &  \text{no}\\
E_8^2 & z^2+x^3+y^5+x^2y^2  & \SL(2,\FF_3) & 8,85 &   \text{no}\\
\end{array}
$$

\vspace{1em}
\centerline{Table \stepcounter{figure}\arabic{figure}: Rational double points of  type $E_n$ in characteristic $p=3$.}
\vspace{1em}

$$
\begin{array}[t]{l|l|l|l|l}
\text{}         &  \text{formal relation}   &  \pi_1^\loc &   l(\O/\shJ),l(\O/\shJ^{[5]}) & \text{$\Theta$ free}\\
\hline&&&&\\[-2ex]
E_6 & z^2+x^3+y^4        & A_4 & 6,173 &  \text{no} \\
\hline&&&&\\[-2ex]
E_7 & z^2+x^3+xy^3        & S_4 & 7,198 & \text{no} \\
\hline&&&&\\[-2ex]
E_8^0 & z^2+x^3+y^5        & 0 & 10,250 &\text{yes}\\
E_8^1 & z^2+x^3+y^5+xy^4  & C_5 & 8,239 &  \text{no}\\
\end{array}
$$

\vspace{1em}
\centerline{Table \stepcounter{figure}\arabic{figure}: Rational double points of  type $E_n$ in characteristic $p=5$.}
\vspace{1em}


\end{document}